\documentclass[reqno,11pt]{amsart}
\usepackage[all]{xy}
\usepackage{graphics,color}
\usepackage{wrapfig}

\usepackage{amssymb,amsmath}

\newtheorem{defn}{Definition}[section]

\newtheorem{obs}{Observation}
\newtheorem{conj}{Conjecture}

\newtheorem{eg}{Example}[section]
\newtheorem{theorem}{Theorem}[section]
\newtheorem{lemma}[theorem]{Lemma}

\newtheorem{cor}[theorem]{Corollary}

\newtheorem{remark}[theorem]{Remark}
\newtheorem{prop}[theorem]{Proposition}
\def\g{{\gamma}}

\newcommand{\BZ}{{\mathbb Z}}

\newcommand{\al}{\alpha}
\newcommand{\G}{{\Gamma}}
\newcommand{\SL}{{SL_2(\BZ)}}

    \let\c\gamma  \let\d\delta
  \let\g\gamma

\def\C{\mathbb C}
\def\G{\mathbf G}

\def\g{\gamma}

\def\Gal{{\rm Gal}}

\def\Z{{\mathbb Z}}

\def\Q{{\mathbb Q}}

\def\G{\Gamma}

\DeclareRobustCommand{\qed}{%
  \ifmmode % if math mode, assume display: omit penalty etc.
  \else \leavevmode\unskip\penalty9999 \hbox{}\nobreak\hfill
  \fi
  \quad\hbox{$\blacksquare$}}

\begin{document}

\title{On  modular forms for some noncongruence
subgroups of $\SL$}
%\shorttitle{On}
\author{Chris A. Kurth}
\address{Department of Mathematics\\Iowa State University\\Ames, IA 50011 \\USA}
\email{kurthc@iastate.edu}
\author{Ling Long}

\address{Department of Mathematics\\Iowa State University\\Ames, IA 50011 \\USA}
\email{linglong@iastate.edu}

%\date{Dec. 1, 2006}
\begin{abstract}In this paper, we  consider modular forms for
finite index subgroups of the modular group whose Fourier
coefficients are algebraic. It is well-known that the Fourier
coefficients of any holomorphic modular form for a congruence
subgroup (with algebraic coefficients) have bounded denominators. It
was observed by Atkin and Swinnerton-Dyer that this is no longer
true for modular forms for noncongruence subgroups and they pointed
out that unbounded denominator property is a clear distinction
between modular forms for noncongruence and congruence modular
forms. It is an \emph{open question}  whether genuine noncongruence
modular forms (with algebraic coefficients) always satisfy the
unbounded denominator property. Here, we give a partial positive
answer to the above open question by constructing special finite
index subgroups of $\SL$ called character groups and discuss the
properties of modular forms for some groups of this kind.
\end{abstract}

\maketitle

\section{Introduction}In  \cite{bass-lazard-serre64} Bass, Lazard,
and  Serre proved that  any  finite index subgroup of $SL_n(\Z)$
with $n>2$ is congruence in the sense that it contains the kernel of
a modulo $q$ homomorphism $SL_n(\Z) \rightarrow SL_n(\Z/q\Z)$ for
some natural number $q$. The story is quite different when $n=2$. In
the 19th century, the existence of noncongruence subgroups of the
modular group $P\SL$ was { in} question until the { affirmative}
results of Fricke \cite{fricke1886} and Pick \cite{Pick1886}. Around
the 1960's more noncongruence subgroups were constructed
\cite[etc]{reiner58, newman65,rankin67}.  In
\cite{Millington691,Millington692}, Millington showed there is a
one-to-one correspondence between finite index subgroups of $P\SL$
and legitimate finite permutation groups described in \cite[Theorem
1]{Millington692}. Using Millington's correspondence, Hsu
\cite{hsu96} gave a concrete method of identifying congruence
subgroups. Indeed noncongruence subgroups predominate congruence
subgroups in $P\SL$ \cite{a-sd,Stothers84}. In \cite{a-sd}, Atkin
and Swinnerton-Dyer initiated a serious investigation on the
properties of noncongruence modular forms using computers. Among
other important observations and theorems, Atkin and Swinnerton-Dyer
pointed out that the Fourier coefficients of certain noncongruence
modular forms have unbounded denominators, which is a clear
distinction between the noncongruence and congruence modular forms.

Let $\G$ be a noncongruence subgroup and $\G^c$ its congruence
closure, namely, the smallest congruence subgroup containing $\G$.
Let (UBD) refer to the following condition { on $\G$:}
\begin{quote}
  \emph{Let $f$ be an arbitrary holomorphic integral weight $k\ge 2$ modular form  for  $\G$  but not for $\G^c$ with algebraic Fourier
  coefficients at infinity. Then the Fourier coefficients of $f$ have
  unbounded denominators.}
\end{quote}
A natural and interesting open question is:
\begin{quote}
\emph{Does every noncongruence   subgroup satisfy the
condition (UBD)}?
\end{quote} To the authors
best knowledge, all data about the known genuine noncongruence
modular forms supports a positive answer to the above question. It
should be { made} clear to the readers that this paper is solely
about the unbounded denominator property of the coefficients of
noncongruence modular forms { and that} Atkin and Swinnerton-Dyer
congruences will not be addressed here, but  will be discussed in a
coming paper by Atkin and the second author \cite{AL07}. For
research in this direction, the readers are referred to the original
paper by Atkin and Swinnerton-Dyer \cite{a-sd}, several important
papers by Scholl \cite[etc]{sch85b, sch88}, and some more recent
papers \cite{lly05,all05, long061}. All groups considered here are
of finite index in $\SL$ unless otherwise specified.

Unlike the approach of Atkin and Swinnerton-Dyer in \cite{a-sd},
which is mainly concerned  with subgroups of $P\SL$ with small
indices, Li, Long, and Yang \cite{lly05} considered modular forms
for a noncongruence subgroup {  defined as follows:}
\begin{defn}
 Given a finite index subgroup $\G^0$ of $P\SL$, a normal subgroup $\G$ of $\G^0$ is called a \emph{character group} of $\G^0$ if
  $\G^0/\G$ is abelian. I.e. there exists  a homomorphism
  \begin{equation}\label{eq:chargp}
    \varphi: \G^0 \rightarrow G,
  \end{equation} where $G$ is a  finite abelian group (written multiplicatively) such that $\G=\ker \varphi$.
\end{defn}From now on such a homomorphism $\varphi$ will be fixed.

{ \begin{defn}\label{defn:types}
  Let $\G$ be the kernel of $\varphi: \G^0 \rightarrow G$ with $G$
abelian. We say $\G$ is a character group of type I if there is a
parabolic element $\gamma \in \G^0$ such that $\varphi(\gamma) \neq
1$. If all parabolic elements $\gamma$ of $\G^0$ have
$\varphi(\gamma) = 1$ we say $\G$ is a character group of type II,
and additionally if all parabolic and elliptic elements of $\G^0$
map to 1 we say $\G$ is of type II(A).
\end{defn}}

For example given any positive prime number $p$, $\G^1(p)$ is a type
II character group of $\G^0(p)$ {  (c.f. Example \ref{eg:G(1)})}.
The main difference between character groups of these two types lies
in their cusp widths and the general concept of \emph{level}
introduced by Wohlfahrt \cite{wohlfahrt64} which extends the
classical level definition for congruence subgroups by Klein. For
any finite index subgroup of $\SL$, its level is the least common
multiple of all cusp widths of the group. {  For index-$n$ type II
character groups $\G$ of $\G^0$, each cusp $c$ of $\G^0$ splits into
$n$ different cusps, say $c_1,\cdots, c_n$ in $\G$. The cusp width
of each $c_i$ is the same as the cusp width of $c$ in $\G^0$.
Therefore, the level of $\G$ remains the same as the level of
$\G^0$. However, this is not true for type I character groups in
general.} {  Also note that any genus 0 subgroup $\G^0$ can be
generated by parabolic and elliptic elements only. Hence there does
not exist any nontrivial type II(A) character group of $\G^0$.}
Later in this paper we will consider those $\G^0$ whose genus is 1
so that results on elliptic curves can be applied. The main result
of this paper is the following theorem which gives a partial
positive answer to the open question above.

\begin{theorem}\label{thm:main}
  Let $\G^0$ be any genus 1  congruence   subgroup whose modular curve has a model defined over a number field. If there
  exists a prime number $p$ such that every index-$p$ type II(A) character group
  of $\G^0$  satisfies the
  condition (UBD), then there exists a positive constant $c$ depending on $\G^0$ such
that for any $X\gg 0$,
\begin{equation}
\#\left  \{\text{Type II(A) char. gp } {\G} \text{ of  } \G^0 \
\big| \ [\G^0:\G]<X,  \G \text{ satisfies (UBD)} \right \}>c\cdot
X^2.
\end{equation}
\end{theorem}  { In comparison, we will shown in Lemma
\ref{lem:countII(A)} that
 \begin{equation}
\lim_{X\rightarrow \infty}\frac{ \# \left \{\text{ II(A) type char.
gp } {\G} \text{ of } \G^0 \ \big| \ [\G^0:\G]<X \right \}}{X^2} =
\frac{\pi^2}{12}.
\end{equation}} The result stated in Theorem \ref{thm:main} can be generalized to other genus cases where the power of $X$ will be changed
accordingly, however we will work with the genus 1 case here. We
will justify in this paper that it is computationally feasible to
verify the conditions for the above theorem. In particular, we
consider the type II(A) character groups of $\G^0(11)$ as Atkin's
calculations in \cite{Atkin67} on modular functions for $\G_0(11)$
(which can be easily turned into modular functions for $\G^0(11)$)
will be very useful to us. We will show the conclusion of Theorem
\ref{thm:main} holds for $\G^0(11)$.

\smallskip

This paper is organized in the following way: in Section
\ref{sec:fieldofmf} we give a general discussion on   modular
functions for character groups where we will present an unbounded
denominator criterion (Lemma \ref{lem:unbounded}). In Section 3, we
will conclude that  {  almost  all nontrivial type II character
groups of the standard congruence subgroups are noncongruence.} In
Section 4, we will construct modular functions for type II(A)
character groups {  in genus 1 } and prove our main result, Theorem
\ref{thm:main}. In Section 5, we will study modular functions for
type II noncongruence character groups of $\G^0(11)$. In the last
section, we will briefly describe type I character groups of
$\G^0(11)$.

The authors are indebted to Prof. A.O.L. Atkin. The authors would
 like to thank him for his communications through which we have
refined our approach. Atkin has provided another proof of the
integrality of the Fourier coefficients of $x$ and $y$ used in
Section 5.  The authors also thank Prof. Wenching Winnie Li, Frits
Beukers, and Siu-Hung Ng for their helpful communications and
constructive suggestions. Prof. Helena Verrill's Fundamental Domain
Drawer was used to generate the fundamental domain in this paper.
The authors are grateful to  the referee whose comments led to a
significantly improved presentation of this   paper.

\section{Fields of modular functions}\label{sec:fieldofmf}

Recall some useful notation and results in \cite{shim1}. An element
$\gamma$ in $P\SL$ is said to be parabolic (resp. elliptic or
hyperbolic) if $|\textrm{tr} \gamma |=2$ (resp. $<2$ or $>2$). {  We
assume $\G^0$ is a congruence  subgroup of $P\SL$ and $\G$  a
character group of $\G$. Denote by $M_k(\G)$ the space of weight $k$
holomorphic modular forms for $\G$. } Let $\C(\G)$ denote the field
consisting of meromorphic modular functions for $\G$ over $\C$. It
is a finite algebraic extension of $\C(\G^0)$. Let $X_{\G}$ denote
the compact modular curve for $\G$. In this paper, we assume
$X_{\G}$ admits an algebraic model. We use standard notation for
some well-known congruence subgroups of $P\SL$ with a given level
$n$. For example
$$\G^0(n)=\left \{ \gamma \in \SL \ \bigg| \ \gamma=\begin{pmatrix}
  *&0\\ *&*
\end{pmatrix} \text{mod } n \right \}\bigg/\pm I_2;$$ likewise we will use $\G^1(n),$ $\G_0(n),$ and $\G(n)$ to denote other standard congruence subgroups.
\begin{lemma}[C.f. Section 2.1 \cite{shim1}]
If $\G$ is a normal subgroup of $\G^0$, then $\C(\G)$ is Galois over
$\C(\G^0)$ with the Galois group $\Gal(\C(\G)/\C(\G^0))$ being
isomorphic to $\G^0/\G$.
\end{lemma}
For any $\gamma \in \G^0$ and $g(z) \in \C(\G)$, $\g$ acts on $g(z)$
via the stroke operator $$f(z)|_{\g}=f(\g z).$$

\begin{lemma}\label{lem:subexten-subgp}
  Normal field extensions of $\C(\G^0)$ which
  are contained in $\C(\G)$  are in one-to-one correspondence with normal subgroups of $\G^0$ containing $\G$.
\end{lemma}
\begin{proof}
By the Galois correspondence,  there is a bijection between normal
intermediate fields between $\C(\G)$ and $\C(\G^0)$   and normal
subgroups of $\G^0/\G$. By one of the isomorphism theorems, normal
subgroups of $\G^0/\G$ are in one-to-one correspondence with normal
subgroups of $\G^0$ which contain $\G$.
\end{proof}

\begin{lemma}\label{lem:intergp2type}
  If $\G$ is normal in $\G^0$ and $\G^0/\G$ is a finite abelian group, then  any
  group $\G'$ sitting between $\G^0$ and $\G$ is normal in $\G^0$.
\end{lemma}
\begin{proof}
  Assume $\G'=\cup_i \d_i \G$. Then for any $\gamma \in \G^0$, $$\gamma
  \G' \gamma^{-1}=\cup_i \gamma \d_i \G \g^{-1}=\cup_i \d_i \g
  \G\g^{-1}=\cup _i \d_i \G=\G'.$$
\end{proof}

\begin{cor}\label{cor:1}
  If $\G$ is any finite index character group of
  $\G^0$, then any intermediate group $\G'$ sitting between $\G^0$ and
  $\G$ is also a character group of $\G^0$. In particular, if $\G$
  is of type I (resp. II, or II(A)), $\G'$ is of the same type.
\end{cor}

In order to discuss the bounded or unbounded denominator property
satisfied by integral weight meromorphic modular forms for $\G$, we
will restrict ourselves to a suitable algebraic number field $K$
rather than $\C$ itself. For simplicity, we use (FS-A) to refer to
the following condition satisfied by a modular form  $f$:
\begin{quote}
 \emph{all the Fourier coefficients of $f$ at infinity
 are in  $K$.}

\end{quote} If $\C(\G)$ is generated over
$\C$ by some transcendental generators $f_1,\cdots, f_n$ where each
$f_i$ satisfies (FS-A), then $K(\G)=K(f_1,\cdots, f_n)$ is a
subfield of $\C(\G)$ consisting of meromorphic modular functions for
$\G$ satisfying (FS-A). Moreover
\begin{equation}
  K(\G)\otimes _K \C=\C(\G).
\end{equation}

Use (FS-B)  to refer to the following condition:
\begin{quote}
 \emph{the Fourier coefficients of $f$ at  infinity  have bounded
 denominators.}
\end{quote}{ Namely, there is an algebraic number
 $C$ such that  the Fourier coefficients of $C\cdot f$ are all
 algebraically integral.}
{In the following discussion, we will assume $K$  is a large number
field. {  Let $$R_K(\G):=\left \{f\in K(\G) \big | f \text{ has at
most a single pole at infinity} \right \}.$$ It is a ring and   its
fraction field is $K(\G)$. }  The following observation is very
crucial to our discussion.
\begin{obs}\label{obs:1}
Elements in $K(\G)$ satisfying (FS-B) are closed under addition,
subtraction, multiplication, and {  hence form a subring of
$R_K(\G)$ which will be denoted by $B_K(\G)$. I.e. $$B_K(\G):=\left
\{ f\in R_K(\G) \big| f \textrm{ satisfies (FS-B)} \right \}.$$
Moreover, for any formal power series $f$ whose leading coefficient
is 1 and Fourier coefficients are all algebraically integral {
(such as a normalized newform)}, $1/f$ also satisfies (FS-B).}
\end{obs}
 {
\begin{lemma}
The ring  $R_K(\G^0)$ is a subset of $B_K(\G)$.
\end{lemma}
\begin{proof}
It is known that that  any congruence cuspform which satisfies
(FS-A) also satisfies (FS-B) (c.f. \cite{shim1}). Given $f\in
R_K(\G^0)$ one can pick a normalized newform $F_1$ for $\G^0$ and a
large positive integer $N$ such that $f\cdot F_1^N$ is a cuspform
for $\G^0$; hence $f\cdot F_1^N$ satisfies (FS-B) and so does $f$ as
$1/F_1^N$ also satisfies (FS-B). We conclude that $B_K(\G)$ contains
$R_K(\G^0)$.  \end{proof}  }

{  Let $$K_B(\G):= \textrm{the fraction field of the ring }
B_K(\G).$$} By the above lemma, $K_B(\G)$ is a field extension over
$K(\G^0)$
\begin{obs}
Let $\G$ be a character group of $\G^0$. By the Fundamental Theorem
of Finite Abelian Groups, it suffices to consider only character
groups such that $\G^0/\G$ is cyclic.
\end{obs}

\begin{lemma}\label{lem:cyclicextn}
  If $\G$ is a character group of $\G^0$ with cyclic quotient of
  order $n$, then there is a modular function $f$ for $\G^0$ such
  that $\C(\G)=\C(\G^0)(\sqrt[n]{f}).$
\end{lemma}
\begin{proof}
Since $\Gal(\C(\G)/\C(\G^0))$ is isomorphic to $\G^0/\G=<a \G>$ for
some  coset $a \G$, there exists a modular function $g \in \C(\G)$
such that $g|_{a}=e^{2\pi i /n}\, g$. Let $f=g^n\in \C(\G^0)$. Then
$\C(\G)$ is a splitting field of  $x^n-f \in\C(\G^0)[x]$.
\end{proof}

 \begin{lemma}\label{lem:2.5}
Assume $\G^0$ is a genus larger than 0 finite index subgroup of
$P\SL$ and $\G$ is a character group of $\G^0$. If
$R_K(\G)=R_K(\G^0)[g]$ for some $g$ satisfying $g^n-f$ where $f\in
R_K(\G^0)$ and $K$ contains the primitive $n$th
   roots,
   then $B_K(\G)=R_K(\G^0)[g^d]$  and $K_B(\G)=K(\G^0)(g^d)$ for some $d|n$.
 \end{lemma}
 \begin{proof}
%By the assumption, $K(\G)$ is the splitting field of the irreducible polynomial $x^n-g$ over $K(\G^0)$. The field $K_B(\G)$  is sitting between $K(\G^0)$ and $K(\G)$. Since  $K(\G)$ is Galois over $K(\G^0))$  with the Galois group $\Gal(K(\G)/K(\G^0))$ being isomorphic to $\Z/n\Z$, by the Galois theory,  the minimal polynomial of $u$ over $K_B(\G)$   is of the form $x^{d}-u^{d}$ for some positive integer $d$ dividing and $K_B(\G)$ is also Galois over $K(\G^0)$ with $u^d$ being a generator.

Given any
 $r=\sum_{i=1}^s a_i g^{n_i}\in B_K(\G)$ with $s\ge 1$, $h_i\neq 0$ for all $i$,  and $1\le n_1<\cdots
 <n_s<n$,  we will show that $g^{n_i}\in B_K(\G^0)$
for all $i\in [1,\cdots, s]$ by mathematical induction.  When $s=1$,
it follows from Lemma \ref{lem:appendix} in the Appendix.

We now assume  $s\ge 2$. Let $H_1$ be a  weight 2 normalized newform
for $\G^0$ so   $1/H_1$ satisfies (FS-A) and (FS-B). One can choose
$H_2\in R_K(\G^0)$ such that  the weight $-2$ meromorphic modular
form $E=H_2/H_1$ and the weight 0 modular function $\displaystyle
\frac{E}{f}\cdot \frac{df}{dz}$ for $\G^0$ both have at most one
pole at infinity. It is easy to see $E$ satisfies both (FS-A) and
(FS-B).
 It is well-known that for any nonconstant
meromorphic weight 0 modular form $F$ for a finite index subgroup of
$\SL$, $\displaystyle \frac{dF}{dz}$
  is a meromorphic weight 2 modular form for the same group.  Thus,  $\mathcal D=E\frac{d}{dz}$ is a
linear mapping on both $R_K(\G^0)$ and $B_K(\G)$ as it preserves the
(FS-B) property. Moreover $\mathcal D g^i=\frac{i}{n}\frac{\mathcal
D f}{f} g^i.$   So $$\mathcal D\left (\sum_{i=1}^s a_i g^{n_i}
\right )=\sum_{i=1}^s b_i g^{n_i}, \text{ where } b_i=\mathcal D
a_i+a_i\cdot \frac{n_i}{n} \frac{\mathcal Df}{f}.$$  If
$\displaystyle \frac{b_i}{a_i}=\frac{b_j}{a_j}$ for some $1\le
i<j\le s$, then $\displaystyle \mathcal D \ln
\frac{a_i}{a_j}=\mathcal D \ln f^{(n_j-n_i)/n}$. This is impossible
as $0< n_j- n_i<n$ and hence $f^{(n_j-n_i)/n}$ is not in $K(\G^0)$.
So $$(b_s- a_s\mathcal D) r =\sum_{i=1}^{s-1} (b_s a_i-a_sb_i)
g^{n_i} \in B_K(\G)$$ and it is not zero. Therefore, by induction,
we know $g^{n_i}\in B_K(\G^0)$ for any $i\in [1,\cdots, s]$. Since
in this case, $g^{n_1+n}\in B_K(\G)$ for any $n\ge 0$ implies
$g^{n_1}\in B_K(\G)$. So, if we  let $d$ be the great common divisor
of $m$ such that $g^m\in B_K(\G^0) $. By Lemma \ref{lem:appendix},
$g^d \in B_K(\G^0) $.  Then
 $K_B(\G)=K(\G^0)(g^d)$.
 \end{proof}

 From now on we will confine ourselves
to character groups of this type unless otherwise specified.
  By the Galois correspondence,
$K_B(\G)\otimes_K \C$, which is an intermediate field between
$\C(\G)$ and $\C(\G^0)$, corresponds to a character group of $\G^0$
containing $\G$. We denote the corresponding character group by
$B_{\G}$, i.e.
$$B_{\G}:= \textrm{the character group such that }
\C(B_{\G})=K_B(\G)\otimes_K \C.$$

 \begin{prop}\label{prop:key}Let  $\G$ be a character group of a
congruence subgroup $\G^0$ of genus larger than 0 such that
$R_K(\G)=R_K(\G^0)[g]$ for some $g$ satisfying $g^n-f$ where $f\in
R_K(\G^0)$. If $B_{\G}=\G^0$, then $\G$ satisfies (UBD).
\end{prop}

\begin{proof}
If $\G$ does not satisfy (UBD), then there exist an integer $j$ and
a  holomorphic modular form $f\in M_j(\G)\setminus M_j(\G^0)$ which
 satisfies (FS-A) and (FS-B).  One can construct a
 function $f'\in R_K(\G)\setminus R_K(\G^0)$
which also satisfies (FS-A) and (FS-B). To see this we first fix
some weight $j'$ normalized nonconstant newform $F_1$ for $\G^0$
with Fourier coefficients in $K$. For a large enough $n$ one can
find a weight $nj'-j$ holomorphic modular form $F_2$ for $\G^0$
satisfying (FS-A) and a $F_3\in R_K(\G^0)$ such that $f'=f\cdot
F_2\cdot F_3/F_1^n\in R_K(\G)$. The functions $f$, $1/F_1^n$, $F_2$,
and $F_3$ satisfy (FS-B), and so does $f'$. While $F_2\cdot
F_3/F_1^n$ is invariant under the action of $\G^0$ and $f$ is not,
so $f'\notin R_K(\G^00$. Therefore $R_K(\G^0)\subsetneq
R_K(\G^0)[f']\subseteq B_K(\G)$. Hence $K(\G^0)(f')$ which is a
nontrivial field extension over $K(\G^0)$ is contained in $K_B(\G)$.
Therefore, $B_{\G}$ is a nontrivial subgroup of $\G^0$; which
contradicts the condition that $B_{\G}=\G^0$.
\end{proof}

{  To end this section, we provide the following  simple criterion
for detecting whether $\sqrt[p]{f}$ satisfies (FS-B) for a given
$f\in K(\G^0)$} and for   $p$ prime.  Note that any nonzero power
series $f$ in $w$ can be easily normalized, up to multiplying a
power of $w$, into the form $f=\sum_{m\ge 0} a_m w^m$ with
$a_{0}\neq 0$.
\begin{lemma}\label{lem:unbounded}
Let $K$ be a number field, $p$ be any prime number and
$$f=a_{0}+\sum_{m\ge 1}a_mw^m, a_m\in K, a_{0} \neq 0$$
such that  for every $m$, $a_m$ is $\wp$-integral for any prime
ideal $\wp$ in $\mathcal O_K$ above $p$. Expand
$\sqrt[p]{f}=\sum_{m\ge 0} b_mw^m$ formally (which is well defined
once we fix a branch for the $p$th root of $a_{0}$).  If there
exists at least one
 $b_m$ such that  ${-\mbox{ord}_{\wp}(b_m/b_0)> \frac{\text{ord}_{\wp}a_{0}}{p}} $, then
$$\limsup_{m\rightarrow \infty} \! -\text{ord}_{\wp}(b_m) \rightarrow \infty.$$
In other words, the sequence $\{b_m\}$ has unbounded denominators.
\end{lemma}
\begin{proof}
  Assume $\{-\text{ord}_{\wp} \left ( b_m/b_0 \right )\}$ has an upper
  bound, say $$\max\{-\text{ord}_{\wp} \left ( b_m/b_0) \right \} =C$$ which is larger than $\frac{\text{ord}_{\wp} a_0}{p}$ by our assumption. Let $m_0$ be the
  smallest positive integer such that $-\text{ord}_{\wp} (b_{m_0}/b_0)=C$.
  Consider the $m_0\cdot p$ coefficient of $(\sum (b_m/b_0) q^m)^p=\sum \left (a_m/a_0 \right )q^m$. Using a standard $p$-adic analysis argument we  derive a contradiction
  $$\text{ord}_{\wp}a_0\ge -\text{ord}_{\wp} (a_{m_0\cdot p}/a_0)=- p \cdot \text{ord}_{\wp}(b_{m_0}/b_0)> \text{ord}_{\wp}
  a_0.
  $$
\end{proof}

\begin{eg}\label{eg:asdunbounded}\cite[4.2.1]{a-sd}. Let
\begin{equation}\label{eq:eta}
\eta(z)=q^{1/24}\prod_{n\ge 1}(1-q^n), \ q=e^{2\pi i z}.
\end{equation} The function
$$\zeta(z)=\left (\frac{ \eta(z)}{\eta(13z)}\right )^2 =q^{-1}(1-2q-q^2+\cdots)\in
\Z[q^{-1},q]]$$   is a Hauptmodul {  of the congruence subgroup
$\G_0(13)$.} Clearly, for any integer $m>2$,
$$\zeta^{1/m}(z)=q^{-1/m}(1-\frac{2}{m}q+\cdots).$$ By the above
lemma, the coefficients of $\zeta^{1/m}(z)$ have unbounded
denominators. The function $\zeta^{1/m}(z)$ is a Hauptmodul of a
genus zero   noncongruence type I character group of $\G^0(13)$.

\end{eg}

\section{Noncongruence character groups of  type II}

{ \begin{defn} A homomorphism $\varphi: \G^0 \rightarrow G$ ($G$ can
be non-abelian) is said to be of type II if it sends all parabolic
elements in $\G^0$ to the identity of $G$.
 \end{defn}}

\begin{lemma}\label{eg:G(1)}Let $p$ be a prime, then $\G^1(p)$ is a type II character group of $\G^0(p)$.
\end{lemma}
\begin{proof}
  Let $\varphi$ be the following homomorphism
  \begin{eqnarray*}
      \varphi: \quad \G^0(p) &\rightarrow& (\Z/p\Z)^{\times} /\pm 1\\
      \pm \begin{pmatrix}   a&b\\c&d
\end{pmatrix}& \mapsto& (a \text{ mod } p) /\pm 1.
  \end{eqnarray*}
  Then $\ker \varphi =\G^1(p).$ Assume $\gamma=\pm \begin{pmatrix}
a&bp\\c&\pm 2-a
\end{pmatrix}$ is a parabolic element in $\G^0(p)$. So $\det
\gamma =1$ implies that $$a\cdot(\pm2-a)=1 \text{ mod } p,$$ hence
$a=\pm 1 \text{ mod } p$. Thus $\varphi(\g)=1$ for every parabolic
element $\g\in \G^0(p)$.
\end{proof}

 \begin{prop}\label{prop:ncc2ndG0}
   Let $\varphi$ be a  homomorphism of   type II from $\G^0(n)$ to another finite group (not necessarily abelian)  whose kernel $\G$ does not contain $\G^1(n)$. Then
    $\G$ is noncongruence.
 \end{prop}
 \begin{proof}
 Let $\G$ be the kernel of such a type II homomorphism $\varphi: \G^0(n) \rightarrow G$  (as we have mentioned in the
 introduction, the level of $\Gamma$ remains $n$). %earlier that the
% level of $\G$ will remain  $n$ {  as each
%  cusp $c$ of $\G^0(n)$ splits into several cusps in $\G$ each of  which has the same cusp width as that of $c$}.
Now we assume that $\G$ is a congruence
 subgroup. Since both $\G^1(n)$ and $\G$ contain $\begin{pmatrix}
   1&0\\1&1
 \end{pmatrix}$, the generator of the stabilizer of 0, so does  $\G \cap \G^1(n)$.
 Therefore the cusp width of $\G \cap \G^1(n)$ at 0 is  1 (while the cusp width of
  $\G(n)$ at 0 is $n$). It follows $[\G \cap
 \G^1(n): \G(n)]\ge n$. On the other hand we have \begin{eqnarray*}
  &n& =[\G^1(n):
 \G(n)]=[\G^1(n):\G \cap \G^1(n)][\G \cap \G^1(n): \G(n)]\\&&\ge [\G^1(n):\G \cap
 \G^1(n)] n.
 \end{eqnarray*} Hence $[\G^1(n):\G \cap
 \G^1(n)]=1$ and $\G^1(n)\subset \G$.
\end{proof}

Using a similar argument, one can obtain that
\begin{prop}\label{prop:ncc2ndG1}
   Let $\varphi$ be any nontrivial  homomorphism of  type II from  either $\G^1(n)$ or $\G(n)$ to another finite group  with kernel $\G$. Then
    $\G$ is noncongruence.
 \end{prop}

\section{Modular functions for type II(A) character groups { in genus 1}}
In this section, we fix $\G^0$ to be a genus one subgroup of $P\SL$,
{ so that} the modular curve $X_{\G^0}$ for $\G^0$ is an elliptic
curve. We will use $O$ to denote the identity (or the origin) of the
elliptic curve. In such a setting, one can apply well-known results
on elliptic curves (c.f. \cite{sil1}).
%{  Note that if $\G^0$ has genus 1 and algebraic,
%$\C(\G^0)=\C(x,y)$ for  $x,y \in R_K(\G^0)$ with  only one pole at
%infinity of order 2 and 3 respectively \cite{kob1}. Hence
%$R_K(\G^0)=K[x,y]$ and we will fix $x$ and $y$ throughout this
%section. }
Let $\text{Div}^0( X_{\G^0})$ denote the set of all
degree zero divisors on the elliptic curve $X_{\G^0}$.
For any $f\in
\C(\G^0)$, let
$$\text{div}(f)=\sum_P n_P(f) (P),$$ where $n_P(f)=\text{ord}_P(f)$ is the order
of vanishing of $f$ at $P$. Then $\text{div}(f)$ is a finite sum
such that $\displaystyle \sum_P n_P=0$. Recall that $D_1, D_2$  are
said to be \emph{equivalent} and are denoted by $D_1 \sim D_2$ if
$D_1-D_2=\text{div}(f)$ for some  $f\in \C(\G^0)$.

There is a natural homomorphism $\pi: \G^0  \rightarrow
H_1(X_{\G^0},\Z)$, the first homology group of $X_{\G^0}$ with
coefficients of $\Z$. For simplicity, we will simply denote
$H_1(X_{\G^0},\Z)$ by $G^0$ (written additively).  Let $\varphi:
\G^0\twoheadrightarrow G$ be any surjective homomorphism. By
Definition \ref{defn:types}, the group $\G=\ker\varphi$ being a type
II(A) character group  is equivalent to the existence of a
surjective homomorphism $\widetilde{\varphi}: H_1(X_{\G^0},\Z)
\twoheadrightarrow G$ such that $\varphi=\widetilde{\varphi} \circ
\pi$. Under our assumption on genus being 1, $G^0$ is a rank-2 free
$\Z$-module, i.e. a rank-2 lattice.

\begin{lemma}\label{lem:subgpsubgroup}
Let $\G^0$ be a  genus 1 finite index subgroup of the modular group,
then type II(A) character groups $\G$ of $ \, \G^0$ are in
one-to-one correspondence with  rank 2 sublattices  $\widetilde{\G}$
of $G^0$.
\end{lemma}
\begin{proof}
  The correspondence is $\G=\ker \varphi \longleftrightarrow \ker
  \widetilde{\varphi}=\widetilde{\G}$.
\end{proof}

$$\xymatrix{1 \ar[r]\ar[rd] & \G\ar[d]_{\pi} \ar[r] & \G^0 \ar[rr]^{\varphi} \ar[d]^{\pi} && G\\
&\widetilde{\G}\ar[r] & G^0\ar[urr]_{\widetilde{\varphi}}&}$$

From now on we will fix the notation that $\widetilde{\G}$ refers to
a finite index subgroup of $G^0$.

\begin{center}
{\sc Diagram 1}
\end{center}
\begin{cor}
Let $p>1$ be a prime number and $\G^0$ be a genus 1 finite index
subgroup of $P\SL$. There are $p+1$ non-isomorphic index-$p$ type
II(A) character groups of $\G^0$.
\end{cor}
\begin{proof}
 The rank-2 lattice $G^0$ has $p+1$ non-isomorphic index-$p$ sublattices.
\end{proof}

Assume $\G$ is a type II(A) character group of $\G^0$ with
$\G^0/\G\cong \Z/n\Z$. By the Hurwitz genus formula (c.f.
\cite{sil1}), $X_{\G}$ also has genus 1. The natural projection map
$X_{\G} \rightarrow X_{\G^0}$ is a degree $n$ isogeny of elliptic
curves.

 By Lemma \ref{lem:cyclicextn}, $\C(\G)=\C(\G^0)(\sqrt[n]{f})$ for some $f\in
\C(\G^0)$. Let $g=\sqrt[n]{f}$.  Consider the divisor
$\text{div}(g)$ in $\text{Div}^0( X_{\G^0})$. By Proposition 3.4 of
\cite{sil1} there exists a unique point $P$ on $X_{\G^0}$ such that
\begin{equation*}
  \text{div}(g)\sim (P)-(O),
\end{equation*} where   $O$ stands for the origin of the elliptic curve $X_{\G^0}$. We fix the cusp infinity to be the origin of $X_{\G^0}$ in the following
discussion. Since $\text{div}(g^n)=\text{div}(f)$ is a principle
divisor, we know $P$ is an $n$-torsion point of $X_{\G_0}$ by Abel's
Theorem.
 Hence one can pick $f=f_P$ to be a
modular function for $\G^0$ satisfying
\begin{equation}\label{eq:fP}
  \text{div}(f_P)=-n\cdot ({O})+n\cdot (P).
\end{equation}

\begin{lemma}\label{lem:fieldag}
  If $X_{\G^0}$ has an algebraic model defined over a number field
  $K_0$ and $\C(\G^0)=K_0(\G^0)\otimes _{K_0}\C$, then for any type II(A) character
  group
  $\G$ of $\G^0$ with $\G^0/\G\cong \Z/n\Z$, $\C(\G)$ is generated over
  $\C(\G^0)$ by a single element $g$ whose Fourier coefficients at
  infinity are in a fixed number field $K$ above $K_0$.
\end{lemma}
\begin{proof}
  Assume $\text{div}(g)\sim (P)-(O)$ and by the above assumption the
  coordinates of $P$ on the elliptic curve $X_{\G^0}$ are algebraic
  numbers. Let $x$ and $y$ be  the local variables of $X_{\G^0}$  which have only one pole at infinity of degree 2 and 3
  respectively \cite{kob1}.
  We may assume the Fourier coefficients of $x$ and $y$ at infinity are in a number field $K_0$. So $R_{K_0}(\G^0)=K_0[x,y]$. Then
  $f_P$ can be expressed as a polynomial in $x$ and $y$ whose
  coefficients can be determined by the local condition that $f_P=F(x,y)$
  has value 0 at $P$ with multiplicity $n$. Hence the
  coefficients of $F(x,y)$ can be determined by solving a system of algebraic
  equations. Then we can pick $K$ to be a large enough number field
  which contains coefficients of $F(x,y)$ and all primitive $n$th roots.
\end{proof}

{  In summary,  if $\G^0$ has genus 1 and $\G$ is a type II(A)
character group of $\G^0$ with cyclic quotient then
\begin{equation}\label{eq:7}
R_K(\G)=R_K(\G^0)[\sqrt[n]{f_P}]
\end{equation} for some $f_P \in R_K(\G^0)$. So
we can apply the results of Lemma \ref{lem:2.5} and Proposition
\ref{prop:key}. By Lemma \ref{lem:2.5}, there is an intermediate
group $B_{\G}$
 between $\G^0$ and $\G$ which corresponds to modular forms for $\G$ satisfying (FS-B).

In the following discussion, for a fixed prime number $p$ we will
construct all non-isomorphic index-$p$ type II(A) character groups
of $\G^0$. }

\begin{lemma}
  Given two functions $g_1$ and $g_2$ in $R_K(\G)$,  assume  $\text{div}(g_1)\sim (P_1)-(O)$ and
  $\text{div}(g_2)\sim (P_2)-(O)$. They generate the same finite field extension over
  $\C(\G^0)$ which corresponds to a type II(A) character group if and only
  if $\langle P_1\rangle=\langle P_2\rangle$ as  finite abelian subgroups of $X_{\G^0}$.
\end{lemma}
\begin{proof}By our previous assumptions, the orders of $P_1$ and
$P_2$ in $X_{\G^0}$ are both $n$.

If $g_1\in \C(\G^0)( g_2)$, then $(P_2)-(O) \in \langle
(P_1)-(O)\rangle$, thus $P_2\in \langle P_1\rangle$. So $\C(\G^0)(
g_1)=\C(\G^0)( g_2)$ implies $\langle P_1\rangle=\langle
P_2\rangle$.

  Conversely, if
  $\langle P_1\rangle=\langle P_2\rangle$ then $P_2=kP_1$, for some integer $k$ such that
  $gcd(k,n)=1.$ So $\text{div}(g_1^k/g_2)=k(P_1)-(kP_1)-(k-1)(O)$
  is principle. Hence $g_2 \in \C(\G^0)(g_1)$. Similarly, $g_1 \in
  \C(\G^0)(g_2).$ Therefore $g_1$ and $g_2$ generate the same field.
\end{proof}

The following proposition follows from the previous discussions.
\begin{prop}\label{prop:fieldextn2type}
 {  Let $\G^0$ be a genus 1 congruence subgroup of $P\SL$}  and $p$ be a prime number. Let $P$ and $Q$ be
  two linearly independent $p$-torsion points of $X_{\G^0}$, then each
  function
  $\sqrt[p]{f}_P, \, \sqrt[p]{f_{Q+iP}}, i=1,\cdots, p$ generates a degree $p$ field
  extension of $\C(\G^0)$ which corresponds to an index-$p$ type II(A) character
  group of $\G^0$. Moreover, any two II(A) character groups obtained
  this way are non-isomorphic.
\end{prop}

{ Next, we  start to estimate the number of  type II(A) character
groups of $\G^0$ satisfying the condition (UBD). } We let
\begin{equation}
  B_{II(A)}=\bigcap _{\G} B_{\G},
\end{equation}where $\G$ runs through all  type II(A) character groups of
$\G^0$. {  The next lemma shows $B_{II(A)}$ is crucial to the
discussion of the (UBD) condition.}
\begin{lemma}\label{lem:BforUBD} {  Let $\G^0$ be a genus 1 congruence subgroup of
$P\SL$}
and $\G$ be an arbitrary   type II(A) noncongruence character group
$\G$ of $\G^0$. If \begin{equation} B_{II(A)} \cdot
  \G=\G^0,
  \end{equation} then $\G$ satisfies the condition (UBD).
\end{lemma}
\begin{proof}The group $B_{\G}$ contains both $B_{II(A)}$ and $\G$.
Since $B_{II(A)}\cdot \G$ is the smallest subgroup of $\G^0$
containing both $B_{II(A)}$ and $\G$, it is contained in $B_{\G}$.
So $$\G^0\supset B_{\G} \supset B_{II(A)}\cdot \G =\G^0$$ and thus
$B_{\G}=\G^0.$ The claim then follows from Proposition
\ref{prop:key}.
\end{proof}

For simplicity, we write  $B$ for $B_{II(A)}$ below.
\begin{lemma}\label{lem:[G:B]finit}
  If there exists a prime number $p$ such that for every index-$p$
  type II(A)
  character group $\G$ of $\G^0$,  $B\cdot \G=\G^0$, then $[\G^0: B]<\infty$.
\end{lemma}
\begin{proof} {  We will stick to previous notation (c.f. Diagram 1 and Lemma \ref{lem:subgpsubgroup}).} Let $\widetilde{B}=\pi(B) \subset G^0.$
  Since $B\cdot \G=\G^0$, we have $\widetilde{B}+ \widetilde{\G}=G^0$. Because $\widetilde{\G}$ is a proper subgroup of $G^0$,
  $\widetilde{B}$ is not trivial.  To achieve the claim it suffices to show that the rank of
  $\widetilde{B}$ is 2. We will rule out the other remaining
  possibility: $\widetilde{B}$ has rank 1.

  Assume $\widetilde{B}$ has rank 1 and, up to picking a new basis
  for $G^0$, we may assume $\widetilde{B}=\langle na \rangle$ for
  some integer $n>0$. Hence $\widetilde{B}+ \langle a, pb\rangle=\langle a, pb\rangle $
  where $\langle a, pb\rangle=\widetilde{\G}$ is an index-$p$ subgroup of $G^0$.
  Let $\G$ be the index-$p$ type II(A) character group of $\G^0$ corresponding
  to $\widetilde{\G}$.
  Correspondingly we have $B \cdot \G=\G$ which contradicts the
  assumption.
\end{proof}

We now fix a set of generators $\{a, b\}$ for the lattice $\G^0$ and
consider  sub-lattices $\widetilde{\G}$ of $G^0$. By a standard
argument using modules over $\Z$, we know each such $\widetilde{\G}$
can be written uniquely as $\langle l a+nb,mb\rangle$ for some
nonnegative integers $l,n,m$ where $0\le n<m$ and $l>0$. Hence
finite subgroups of $G^0$ are in one-to-one correspondence with
triples $(l,n,m)$ of nonnegative integers satisfying $0\le n<m$ and
$l>0$.

 We  first give an estimation for the number of type
II(A) character groups  of $\G^0$  as follows.
\begin{lemma}\label{lem:countII(A)}
    \begin{equation}
\lim_{X\rightarrow \infty}\frac{ \# \left \{\text{ II(A) type char.
gp } {\G} \text{ of } \G^0 \ | \ [\G^0:\G]<X \right \}}{X^2} =
\frac{\pi^2}{12}.
\end{equation}
\end{lemma}

\begin{proof}Let $X$ be a large positive integer. Let $$S(X):=\# \left \{\text{ II(A) type char.
gp } {\G} \text{ of } \G^0 \ \big| \ [\G^0:\G]<X \right \}.$$
Computing
 $S(X)$ boils to counting  the number  of triples of non-negative
integers $(l,m,n)$ such $l\cdot m<X, 0\le n<m$. Let $l$ be a fixed
positive integer smaller than $X$, then $m$ can be any positive
integer no bigger than $\lceil \frac{X}{l} \rceil$ and $n$ can be
any nonnegative integer smaller than $m$, so there are altogether
 \begin{equation}
  S(X)=\sum _{l=1}^{X-1} \sum_{m=1}^{\lceil \frac{X}{l} \rceil-1} m=\sum
  _{l=1}^{X-1}\frac{\lceil \frac{X}{l} \rceil (\lceil \frac{X}{l}
  \rceil-1)}{2}.
\end{equation}
Hence \begin{equation}\lim_{X\rightarrow \infty}\frac{ S(X)}{X^2} =
\frac{1}{2} \zeta(2)=\frac{\pi^2}{12},
\end{equation} where $\zeta(2)$ is the value of the classical
Riemann-Zeta function at 2.
\end{proof}

\begin{lemma}
Assume $\widetilde{B}=\langle s a+ub,vb\rangle$ with such a triple
$(s,u,v)$. Then $\widetilde{\G}+\widetilde{B}=G^0$ if and only if
\begin{equation}\label{eq:mnl}
gcd(s,l)=1=gcd(v,m,sn-ul).
\end{equation}
\end{lemma}
\begin{proof} C.f. Theorem 3.9 \cite{jacobson1}.
\end{proof}

Next we give an estimation for the number of  type II(A) character
groups $\G$ of $\G^0$ satisfying $\G \cdot B=\G^0$. By Lemma
\ref{lem:BforUBD}, those $\G$ satisfy the condition (UBD).

\begin{lemma}\label{lem:4.10}Assume  $\G^0$ is a genus 1 congruence subgroup and $[\G^0:B]<\infty$. There exists a positive constant $c$ depending on $\G^0$ such that
when $X\gg0$
\begin{equation}
 \lim_{X \rightarrow \infty} \frac{\# \left \{\text{ II(A) type char. gp } {\G} \text{ of  } \G^0  \ |
\ [\G^0:\G]<X,  \ \G\cdot B=\G^0 \right \}}{S(X)} >c.
\end{equation}
\end{lemma}
\begin{proof}
  It is equivalent to consider finite index subgroups $\widetilde{\G}$ of $G^0$ such that $\widetilde{\G}+\widetilde{B}=G^0$.
  By the previous lemma, it boils down to counting the number of triples $(l,n,m)$ satisfying  $0\le n<m$, $l>0$, and \eqref{eq:mnl}.
 Now we count the number of triples such that $0\le n<m$, $l>0$, $m\cdot l<X$ and
 $gcd(s,l)=1=gcd(v,m)$. We further assume $l=1$ and $X/2<m<{X}$. The problem has been reduced to
  counting the number of couples $(m,n)$ such that $$
  X/2<m<{X}, \quad
  (s,m)=1,\quad  0\le n\le m.$$ Between $X/2$ and $X$ there are about $\frac{\phi(s)}{2s}X$ integers coprime to $s$, where
  $\phi(s)$ is the Euler number of $s$. So there exists a constant $c_1>0$ depending on $s$  such that
there are at least $c_1\cdot X$  positive integers
  within $(X/2,X)$
  satisfying $(s,m)=1$. Therefore there are at least
  $\frac{c_1}{2}\cdot X^2$ triples of nonnegative integers  $(l,n,m)$
  satisfying the conditions above. The statement of the Lemma then
  follows from
  \ref{lem:countII(A)}.

\end{proof}
\begin{proof}[Proof of Theorem \ref{thm:main}]

{ Let $\G^0$ be a genus 1 congruence subgroup  and $p$ be a prime
number. Each index-$p$ type II(A) character group $\G$ of $\G^0$
satisfies $R_K(\G)=R_K(\G^0)[\sqrt[p]{f_P}]$ for some function
$f_P\in R_K(\G^0)$ (c.f. \eqref{eq:fP} and \eqref{eq:7}). If
$\sqrt[p]{f_P}$ does not satisfy (FS-B), then by Lemma
\ref{lem:2.5},  $R_K(\G)=R_K(\G^0)$ and  $B_{\G}=\G^0$. It follows
from Proposition \ref{prop:key} that $\G$ satisfies (UBD). If all
$p+1$ non-isomorphic index-$p$ type II(A) character groups of $\G^0$
satisfy (UBD), then Lemma \ref{lem:[G:B]finit} implies
$B_{II(A)}=\bigcap _{\G} B_{\G}$ is a finite index subgroup of
$\G^0$. Our main result, Theorem \ref{thm:main}, then follows from
Lemma \ref{lem:4.10}. }\end{proof}

\section{Character groups of $\G^0(11)$ of type II} In this chapter, we show that it is computationally feasible to verify the
conditions of Theorem \ref{thm:main}  by working with type II(A)
character groups of $\G^0(11)$. We choose $\G^0(11)$ as the
integrality of the Fourier coefficients of two basic modular
functions for $\G^0(11)$ is known due to a result of Atkin
\cite{Atkin67}.
\subsection{The group $\G^0(11)$}

\begin{figure}
\begin{center}
\scalebox
{.5} % h_length
{
\includegraphics*{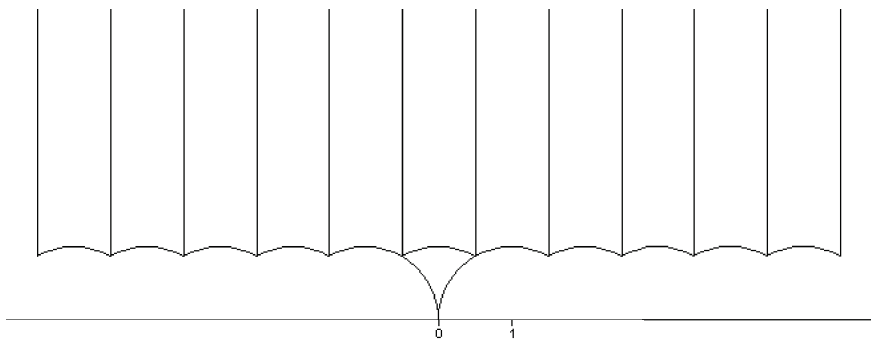}
}
\end{center} \caption{Fundamental domain for
$\G^0(11)$ }
\end{figure}

{  The group $\G^0(11)$ has genus 1 and  is torsion free. A set of
generators of $\G^0(11)$ can be chosen as  two parabolic elements
$\g_{\infty}$ and $\g_0$ and two hyperbolic elements  $A_1$ and
$B_1$ subject to only one relation:
$$\g_{\infty}\g_0A_1B_1A_1^{-1}B_1^{-1}=I_2.$$ Since $\G^0(11)$ does
not have any elliptic elements, every type II character group of
$\G^0(11)$ is automatically of type II(A).}

The modular curve  $X_{\G^0(11)}$ has an equation (c.f.
\cite{cremona-book})
\begin{equation}\label{eq:1}
 X_{\G^0(11)}: y^2 + y = x^3 - x^2 - 10x - 20.
\end{equation}

It is known \cite[III, Prop. 19]{kob1} that the unique (up to a
scalar) holomorphic differential 1-form on $X_{\G^0(11)}$ is
\begin{equation}\label{eq:dx/y}
\frac{dx(z)}{2y(z)+1}=c \cdot \eta(z)^2 \eta(z/11)^2 dz.
\end{equation} %where
%\begin{equation}
%  \eta(z)=q^{1/24}\prod_{n\ge 1}(1-q^n),\quad q=e^{2\pi i z}.
%\end{equation}
One can easily determine that $c=1$.
From \eqref{eq:1} and \eqref{eq:dx/y}, we have
$$x(z):=w^{-2}+2 w^{-1}+4+5 w+8 w^2+ w^3+7 w^4-11 w^5+10 w^6-12 w^7-
18 w^8+\cdots,$$
$$y(z):=w^{-3}+3 w^{-2}+7 w^{-1}+12+17 w+26 w^2+19 w^3+37 w^4-15 w^5-
16 w^6-67 w^7+\cdots,$$ where $w=e^{2\pi i z/11}$.

\begin{lemma}
  The Fourier coefficients of $x$ and $y$ at infinity are in $\Z$.
\end{lemma}
\begin{proof}
In \cite{Atkin67}, Atkin gave a system of modular functions
$\{G_n(z)\}$ on $X_{\G_0(11)}$ (and $\{G_n(z/11)\}$ on
$X_{\G^0(11)}$) with integral coefficients and only one pole of
order $n$ at infinity and a zero at $z=0$ with maximal
multiplicities. According to Lemma 4 of \cite{Atkin67}, one can
obtain that
$$\begin{array}{lll}
x&=&G_2(z/11) + 16\,,\\
y&=&G_3(z/11) + 6\,G_2(z/11) + 60\,.\\
\end{array}
$$
Therefore $x$ and $y$ have integral coefficients.
\end{proof}

\begin{cor}
  The $x$ and $y$ coordinates for $z=0$ are 16 and 60.
\end{cor}

\subsection{Index-2  type II character groups of $\G^0(11)$}
By direct verification, we can describe the index-2 type II
character groups of $\G^0(11)$ as follows.
\begin{lemma}
There are three distinct surjective homomorphisms $\varphi:
\G^0(11)\rightarrow \{\pm 1\}$ of  type II.
\begin{itemize}
  \item If $\varphi(A_1)=-1, \varphi(B_1)=1$, then $\ker \varphi$ is generated
  by $$\{A_1^2, B_1, \gamma_0,\gamma_{\infty}, A_1 \gamma_0 A_1^{-1},
  A_1
  \gamma_{\infty} A_1^{-1}\}$$ subject to the relation
  $$ (\gamma_0 \gamma_{\infty} )(A_1 \gamma_0 \gamma_{\infty} A_1^{-1} )(A_1^2
  B_1
  A_1^{-2} B_1^{-1})=I_2.$$\
  \item If $\varphi(A_1)=1, \varphi(B_1)=-1$, then $\ker \varphi$ is generated
  by $$\{A_1, B_1^2, \gamma_0,\gamma_{\infty}, B_1 \gamma_0 B_1^{-1},
  B_1
  \gamma_{\infty} B_1^{-1}\}$$ subject to the relation
  $$ (B_1 \gamma_0 \gamma_{\infty} B_1^{-1})( \gamma_0 \gamma_{\infty})(
  A_1
  B_1^2
  A_1^{-1} B_1^{-2})=I_2.$$
  \item If $\varphi(A_1)=-1, \varphi(B_1)=-1$, then $\ker \varphi$ is generated
  by $$\{A_1^2, B_1A_1, \gamma_0,\gamma_{\infty}, A_1 \gamma_0 A_1^{-1},
  A_1
  \gamma_{\infty} A_1^{-1}\}$$ subject to the relation
  $$ (\gamma_0 \gamma_{\infty})( A_1 \gamma_0 \gamma_{\infty} A_1^{-1})( (A_1^2)
  (BA_1)
  (A_1^2)^{-1} (BA_1)^{-1})=I_2.$$
\end{itemize}
\end{lemma}

\noindent Making the change of variables
$$x:=\frac{1}{\sqrt[3]{4}}X+\frac{1}{3},\quad y:=\frac{1}{2}Y-\frac{1}{2},$$ equation \eqref{eq:1} will be changed to
$$Y^2=X^3-\frac{31}{3}\sqrt[3]{2^4}X-\frac{2501}{27}.$$
Let $\al_i, {1\le i \le 3}$ be the three roots of
$X^3-\frac{31}{3}\sqrt[3]{2^4}X-\frac{2501}{27}$. Then $X-\al_i$ is
a meromorphic function on $X_{\G^0(11)}$  and
$$\text{div}(X-\al_i)=2(P_i)-2(O),$$ where the $(X,Y)$ coordinates for  $P_i$ are $(\al_i,0)$. Let
$\beta_i=\sqrt[3]{4}/3+\al_i$. Then
$$f_{P_i}=X-\al_i=\sqrt[3]{4}x-\beta_i=\sqrt[3]{4}w^{-2}+2\sqrt[3]{4}
w^{-1} +\cdots.$$ Formally
$$\sqrt{f_{P_i}}=\sum b_n w^n=\sqrt[3]{2} \left (w^{-1}+1+(\frac{3}{2}-\frac{\beta_i}{2})w+(1-\frac{\beta_i}{2})w^2+\cdots \right ).$$
Either  $-\text{ord}_{\wp}( b_1/b_{-1})$ or $-\text{ord}_{\wp}
(b_2/b_{-1})$ equals $\text{ord}_{\wp}2$ which is bigger than
$\text{ord}_{\wp} (\frac{\sqrt[3]{4}}{2})$ for some place $\wp$
above $2$.
 Hence by Lemma \ref{lem:unbounded},  {  we derive that
 the coefficients of
 $\sqrt{f_{P_i}}$ have unbounded denominators for $i=1,2,3$. By
 Proposition \ref{prop:key}, we have}
\begin{theorem}
  Theorem \ref{thm:main} holds for $\G^0(11)$.
\end{theorem}
\begin{proof}
  Let $\G$ be any one of the index-2 type II(A) character groups of $\G^0(11)$. The field $\C(\G)$ is an extension of $\C(\G^0(11))$
  generated by an element whose Fourier expansions have unbounded
  denominators. Thus $\G$ satisfies the condition $\G \cdot B_{II(A)}=\G^0(11)$  and hence
  (UBD). So $\G^0(11)$ satisfies the conditions of Theorem
  \ref{thm:main}.
\end{proof}

%\begin{theorem}
%There exists a constant $c$ such that the conclusion of Theorem
%\ref{thm:main} holds for $\G^0(11)$.
%\end{theorem}

\subsection{Index-5 character groups of $\G^0(11)$ of  type II}
In this subsection, we consider index-5 type II character groups of
$\G^0(11)$. Note that $\G^1(11)$ is one of these groups and the
Mordell-Weil group of $X_{\G^0(11)}$ over $\Q$ is isomorphic to
$\Z/5\Z$.

Let $P=[5,5]$ and $Q=[-\frac{1}{2}+\frac{11}{10}\sqrt{5},
-\frac{1}{2}+\frac{11}{10}\sqrt{-25-2\sqrt{5}}]$. They generate all
5-torsion points of $ X_{\G^0(11)}.$ (In particular, $3P=[16,60]$
corresponds to the point $z=0$ on $X_{\G^0(11)}$.) The
$x$-coordinates of $Q+iP, \, {1\le i \le 4}$ are the roots of the
monic polynomial $x^4+x^3+11x^2+41x+101$ and hence are all integral.
Similarly, the $y$-coordinates of $Q+iP, \, {1\le i \le 4}$ are also
integral.

By an explicit calculation using Maple, it is found that (c.f
\eqref{eq:fP})
$$\begin{array}{lll}
f_P&=&xy-4x^2+30x-4y-55.
\end{array}$$ Hence the Fourier coefficients of the expansion of $f_P$ in terms of
$w$ are all 5-integral. When $i=1,2,3,4$, the coefficients of
$f_{Q+iP}$, as a polynomial in $x$ and $y$, are in a larger number
field $K$. However, by checking the denominators, one concludes that
the Fourier coefficients of the expansion of $f_{Q+iP}$ in terms of
$w$ are all $\wp$-integral for any prime $\wp$ in $K$ above 5.

\begin{cor}
  The $w$-expansions of all above functions are $\wp$-integral for any prime $\wp$ above 5.
\end{cor}
More explicitly, we have
$$\begin{array}{lll}
f_P&=&w^{-5}+ w^{-4}-3 w^{-3}+13 w^{-2}+20 w^{-1}-23+\cdots\,,\\
f_{Q+P}&=&w^{-5}+
w^{-4}+\frac{23+\sqrt{5}+i(3+\sqrt{5})\sqrt{25+2\sqrt{5}}}{4}w^{-3}+\cdots\,,\\
f_{Q+2P}&=&w^{-5}+
w^{-4}+\frac{99-33\sqrt{5}+i(23+3\sqrt{5})\sqrt{25+2\sqrt{5}}}{44}w^{-3}+\cdots\,,\\
f_{Q+3P}&=&w^{-5}+
w^{-4}+\frac{99-33\sqrt{5}-i(23+3\sqrt{5})\sqrt{25+2\sqrt{5}}}{44}w^{-3}+\cdots\,,\\
f_{Q+4P}&=&w^{-5}+
w^{-4}+\frac{23+\sqrt{5}-i(3+\sqrt{5})\sqrt{25+2\sqrt{5}}}{4}w^{-3}+\cdots\,.
\end{array}
$$

Hence one can apply Lemma \ref{lem:unbounded} and conclude that the
coefficients of $\sqrt[5]{f_P(w)}$ and $\sqrt[5]{f_{Q+iP}(w)}, \,
1\le i\le 4$ have unbounded denominators. {  Recall also these
functions correspond to 5 non-isomorphic index-5  type II(A)
character groups of $\G^0(11)$ (c.f. Proposition
\ref{prop:fieldextn2type}). The group $\G^1(11)$ is another index-5
type II(A) character group of $\G^0(11)$ corresponding to $f_Q$. }

\begin{theorem}\label{thm:index5}
  There are 6 index-5  type II(A) character groups   in $\G^0(11)$. Among them, one is $\G^1(11)$ and
  the other 5 are noncongruence. Moreover every one of these noncongruence
  subgroups satisfies the condition (UBD).
\end{theorem}

{  The following conjecture is equivalent to the condition (UBD)
being held by all type II(A) character groups of $\G^0(11)$.}
\begin{conj}
  For $\G^0(11)$, $B_{II(A)}=\G^1(11)$.
\end{conj}

\section{Character groups of $\G^0(11)$ of type I }

By Lemma 3 of \cite{Atkin67}, the following  modular function for
$\G^0(11)$ \begin{equation} G_5(z)=\left
(\frac{\eta(z/11)}{\eta(z)}\right )^{12}=w^{-5}\left (1-12 w+54
w^2-88 w^3-99 w^4+\cdots \right )\end{equation} satisfies
$$\text{div}(G_5)=5(P_0)-5(O),$$ where $P_0=[16,60]$ is a 5-torsion point of $X_{\G^0(11)}$. Since $P_0$ corresponds
to the point $z=0$ as we have mentioned earlier, $O$ and $P_0$ are
the two cusps of $\G^0(11)$.

Let $\G$ be a type I character group of $\G^0(11)$ with
$\G^0(11)/\G\cong \Z/n\Z$. Then $\C(\G)=\C(x,y)(\sqrt[n]{f})$ for
some  $f\in \C(\G^0(11))=\C(x,y)$. Consider
$$\text{div}(f)=n_{P_0}(f)(P_0)-n_O(f)(O)+\sum_{i=1}^l n_i(f)(P_i), $$
where each $P_i$ is different from $O$ or $P_0$. Since the covering
map $X_{\G} \rightarrow X_{\G^0(11)}$ only ramifies at the cusps,
$n|n_i \ {1\le i\le l}.$ As $P_0$ is a 5-torsion point, we have
$$\text{div}(f^5)=\text{div} (G_5^{n_{P_0}(f)})+ \left ( \sum_{i=1}^l 5n_i(f)(P_i)+(5n_{P_0}(f)-5n_O(f)) (O) \right ).$$ If $5 \nmid n$, then $\sqrt[n]{f}$ and $\sqrt[n]{f^5}$ will generate
the same field extension over $\C(\G^0(11))$. Hence we may assume
$f=(G_5)^{ n_{0}} \cdot f'$ for some integer $n_0$ and $f' \in
\C(\G^0(11))$ such that $\C(\G^0(11))(\sqrt[n]{f'})$ corresponds to
an index-$5n$ type II character group of $\G^0(11)$.

In particular we consider the case when $f'=1$ and $5\nmid n$. In
this case we may assume $n_0=1$. Let $\G_n$ denote the character
group whose field of modular functions is generated by
$\sqrt[n]{G_5}$ over $\C(\G^0(11))$. By using the genus formula
again, $\G_n$ has genus $n$. When $n\neq 1,2,3,4,6,12$ and $5\nmid
n$, the coefficients of $\sqrt[n]{G_5}$ have unbounded denominators
and hence $\G_n$ is a type I character group of $\G^0(11)$ which
satisfies the condition (UBD). Moreover, $\sqrt[5]{G_5}$ corresponds
to an index-5 type II character group.

\begin{remark}
{ When $n=2,3,4,6,12$, $\sqrt[n]{G_5}$ is an eta quotient. It is
clear that it is a congruence modular function. Hence $\G_n$ is
congruence when $n$ is divisible by 12.}  By the classification of
Cummins and Pauli on congruence   subgroups with genus no larger
than 24 \cite{cummins-pauli03}, we know when $n=2,3,4,6,12$, the
corresponding character groups are $22A^2, 33A^3, 44A^4, 66B^6,
132A^{12}$ in the notation of Cummins and Pauli.
\end{remark}

\section*{Appendix}

\begin{lemma}\label{lem:appendix}
If $h(w)$ and $g(w)$ are formal power series in $w$ with Fourier
 coefficients in a number field $K$ such that  for some integers $1\le n_1<n$,
$h(w)$, $g^n(w)$, and $h(w)\cdot g^{n_1}(w)$ satisfy (FS-B), then
$g^{n_1}(w)$ also satisfies (FS-B).
\end{lemma}
\begin{proof}
We may assume $g(w)$ is not a monomial in $w$ and $n_1$
 and $n$ are relatively prime to each other. Otherwise, we can use
 $g^{d}$ for $g$ and $n/d$ (resp. $n_1/d$) for $n$ (resp. $n_1$).

Let $f(w)=g^n(w).$ We assume $h(w)=\sum_{m\ge 1} c_mw^m$ with all
algebraically integral coefficients and satisfying $c_1\neq 0$ and
$\gcd (c_1,c_2,\cdots)=1$. Using the factorization of $c_1$ in $K$,
there exists $m_0$ and a linear combination $\sum_{i=1}^{m_0}\beta_i
c_i=1$ for some algebraic integers $\beta_i$ in $K$. We note that
$hg^{n_1+n}$ satisfies (FS-B), and so does
$$\frac{d }{dw}hg^{n_1+n} =\frac{dh}{dw}\cdot g^{n_1+n}+
\frac{n_1+n}{n}\cdot \frac{df}{dw} \cdot (hg^{n_1}).$$  It follows
$\displaystyle \frac{dh}{dw} \cdot g^{n_1+n}$ also satisfies (FS-B).
Let $h^{[1]}=h$. Then
$$h^{[2]}= w\frac{dh}{dw}-h = \sum_{m\ge 2 } c_m(m-1) w^m$$ has algebraically integral coefficients and its product with $g^{n_1+n}$ satisfies (FS-B).
We iterate the above process by replacing $h^{[1]}$ by $h^{[2]}$ and
remark that the coefficients of $\displaystyle \frac{d h^{[2]}}{dw}$
have a common factor of 2. So $$h^{[3]}= \frac{w}{2} \frac{d
h^{[2]}}{dw}-h^{[1]}= \sum_{m\ge 3 } c_m\frac{(m-2)(m-1)}{2} w^m.$$
By  induction for any positive integer $k$, the product of
$$h^{[k]}= c_kw^k+ \sum_{m> k } c_m\binom{m-1}{k-1} w^m$$ and
$g^{n_1+n}$ satisfies (FS-B).   Then $H(w)=\sum_{i=1}^{m_0} \beta_i
 h^{[i]}w^{1-i}$ is a power series starting with
$w$ and having all coefficients algebraically integral.  So
$H(w)^{-1}$ satisfies (FS-B). Since $H(w)\cdot g^{n_1+n}(w)$
satisfies (FS-B), so does
 $g^{n_1+n}(w)$.

Since  $\gcd(n,n_1)=1$,  there exists a positive integer $k$ such
that $n_1\cdot k = nl+1$ for some  $l\in\Z$. Hence
$(g^{n_1+n})^k=g^{n(k+l)}g$ satisfies (FS-B). By
 the previous paragraph, $g^{n+1}$ satisfies (FS-B) and so does $g^{n'}$ for any integer
 $n'\ge  n$. We  assume $n_1=1$ in the following discussion.
Up to multiplying a scalar,  $\displaystyle
f(w)=a_{k_0}w^{k_0}+a_{k_0+1}w^{k_0+1}+\cdots$ where $a_{k_0}\neq
0$, $a_i$ are all algebraic integers, and $\gcd
(a_{k_0},a_{k_0+1},\cdots )= 1$. For any non-monomial  formal power
series like $f$ with coefficients in  $K$ and satisfying (FS-B), we
let $m(f)=a_0$. It is a function measuring the least common
denominators of the coefficients of $f$ when $f$ is normalized to
have leading coefficient 1. The function $m(f)$ is well-defined up
to multiplying a unit in  $K$ and $m(f_1f_2)=m(f_1)m(f_2)$ for any
eligible $f_1$ and $f_2$. In the discussion below, we assume
$k_0=0$. Then
$$\frac{f(a_0w)}{a_0}=1+a_1w+\sum_{i\ge 2} a_ia_0^{i-1}w^i$$ has
algebraically integral coefficients and so does $\displaystyle \left
(\frac{f(a_0w)}{a_0}\right )^{-1}$. By our assumptions, since
$g^{n+1}(a_0w)=f(a_0w)\cdot g(a_0w)$ satisfies (FS-B), so does
$g(a_0w)$. Write
$$\frac{g(a_0w)}{\sqrt[n]{a_0}}=1+b_1w+b_2w^2+\cdots.$$ So
$$\frac{g(w)}{\sqrt[n]{a_0}}=1+\frac{b_1}{a_0}w+\frac{b_2}{a_0^2}w^2+\cdots.$$
Suppose $j$ is the least positive integer such that $b_j$  is not
algebraically integral, then for any integer $q>n$ which is
relatively prime to $a_0$ the $j$th coefficient of $g^q(w)$ is not
algebraically integral either.  So $m(g^q)$ is not a multiple of
$(a_0)^{q/n}$ by a unit, which contradicts
$m((g^q)^n)=m((g^n)^q)=a_0^q$.   Hence the coefficients of
$g(w)/\sqrt[n]{a_0}$ are all algebraically integral.

\end{proof}

%\bibliographystyle{amsalpha}
%\bibliography{longbibl}

\providecommand{\bysame}{\leavevmode\hbox
to3em{\hrulefill}\thinspace}
\providecommand{\MR}{\relax\ifhmode\unskip\space\fi MR }
\providecommand{\MRhref}[2]{%
  \href{http://www.ams.org/mathscinet-getitem?mr=#1}{#2}
} \providecommand{\href}[2]{#2}

\end{document}